# The Diagonal Method and Hypercomputation

Toby Ord and Tien D. Kieu

The diagonal method is often used to show that Turing machines cannot solve their own halting problem. There have been several recent attempts to show that this method also exposes either contradiction or arbitrariness in other theoretical models of computation that claim to be able to solve the halting problem for Turing machines. We show that such arguments are flawed — a contradiction only occurs if a type of machine can compute *its own* diagonal function. We then demonstrate why such a situation does not occur for the methods of hypercomputation under attack and why it is unlikely to occur in any other serious methods.

## 1   Introduction

In recent papers by Cotogno [2003] and Svozil [1998], it is claimed that certain theoretical methods for hypercomputation must fail on strictly logical grounds. Their arguments are similar, using a modified form of the diagonal construction due to Turing [1936] to show that the ability to solve the halting problem leads to contradiction. This argument is very familiar, with versions of it appearing in many introductory texts to show that Turing machines cannot solve their own halting problem. To see exactly what this argument tells us about hypercomputation, let us review it in the form given by Cotogno.

Let $\psi_1, \psi_2, \psi_3\ldots$ be an enumeration of the computable (partial)[1] functions and define the halting function:

$f(x, y)$:      **if** $\psi_x(y)$ diverges **then** 0 **otherwise** 1

and the diagonal function:

$g(x)$:      **if** $f(x, x) = 0$ **then** 0 **otherwise** diverge

If $g(x)$ is computable, it must have an index $i$ such that $\psi_i(x) = g(x)$. This implies that $\psi_i(i) = 0 \Leftrightarrow g(i) = 0 \Leftrightarrow f(i, i) = 0 \Leftrightarrow \psi_i(i)$ diverges which is a clear contradiction. Thus, $g(x)$ cannot be computable.

This argument holds for all definitions of computability. It simply says that for any sequence $\{\psi_i\}$ of functions from N to N, one can construct a function $g(x)$ that is different to each $\psi_i(x)$ and thus cannot be a member of that sequence. So far, the notion of computability has played no part.

When the argument is applied to Turing machines, we let $\{\psi_i\}$ be a sequence of all Turing machine computable functions and see immediately that $g(x)$ cannot be computable by a Turing machine. We then typically prove that all steps in the construction of $g(x)$ can be carried out by Turing machines except for computing $f(x,y)$. Since the existence of a Turing machine that computes $f(x,y)$ will imply the existence of a Turing machine for computing $g(x)$, $f(x,y)$ must also be uncomputable by Turing machines. Since these functions $g(x)$ and $f(x,y)$ were constructed relative to a sequence of all Turing machines, we call them the diagonal function for Turing machines and the halting function for Turing machines respectively.

---

[1] For the sake of clarity, we will often just say 'function' when we refer to both partial functions and total functions. If a claim applies only to total functions, this will be made clear.

So how does this argument affect hypercomputation? If we let $\{\psi_i\}$ be a sequence containing all the functions computable by a given type of hypermachine, we can see that it cannot compute its own diagonal function. Thus, if one could show that some proposed type of hypermachine is so powerful that it *could* actually compute its own diagonal function, it must be logically inconsistent.

This is a strong argument against such models of hypercomputation. As Svozil and Cotogno show, this would force one to dramatically revise the model. Svozil suggests that one could escape the contradiction of applying $g(x)$ to its own index, by using qubits instead of classical bits to give a value that is the quantum superposition of a state representing convergence and one representing divergence. Similarly, Cotogno states that one must choose between accepting either output or accepting one of them by convention. However, we shall show that this dilemma for prospective hypermachines dissolves when we consider whether any proposed hypermachines actually claim to be able to compute their own diagonal functions.

Most proposed models of hypercomputation can compute the halting function for Turing machines. While this is not essential to the concept of hypercomputation (it suffices to compute *any* non-recursive function), it is a most prevalent ability. Most of these hypermachines can also go a step further and compute $g(x)$ for Turing machines. Does this lead to logical inconsistency?

No. Firstly, we must stress that $g(x)$ for Turing machines is a perfectly well defined (if not constructive) function from N to N. Along with other non-recursive objects such as $\Omega$ (Chaitin [1987]), $g(x)$ can be easily defined in mathematical theories such as ZFC. Using ZFC, one can also prove many things about $g(x)$, such as that it exists (in the mathematical sense), that it equals 0 for infinitely many values of $x$ and that no Turing machine can compute it, however, ZFC cannot prove which values of $x$ make $g(x)$ converge. This is not a deficiency in our mathematical apparatus for defining such functions, but just in our present ability to assign truth-values to certain statements expressible in this mathematical language.

The behaviour of a hypermachine that could compute $g(x)$ for Turing machines is therefore perfectly well defined: this machine would simply halt on input $x$ (and return 0), if and only if the $x$th Turing machine does not halt on input $x$. This hypermachine need not lead to contradiction by being given its own index as input, because it is only defined to work on Turing machine indices and is not a Turing machine. By running this machine for $g(x)$ and a Turing machine for $\psi_x(x)$ until one of them halts, we could find a 'proof' as to whether $g(x)$ diverges or converges for any given value of $x$. This would extend the set of truths provable in ZFC, but not contradict them.

The fact that $g(x)$ for Turing machines is well defined and computable by a sufficiently powerful machine can be made more clear through analogy with a more simple diagonal function. Let us consider $\{\psi_i\}$ to be a sequence of all primitive recursive functions. Since all primitive recursive functions are total, their $g(x)$ would simply be divergent on all input and is thus trivially not primitive recursive. As this provides no particular insight, we consider a diagonal function that is total:

$h(x)$: $\quad \psi_x(x) + 1$

This new function $h(x)$ also diagonalises out of the primitive recursive functions and thus cannot be primitive recursive. Note that since primitive recursive functions always halt, this cannot be due to a halting problem, but is simply due to the fact that $\psi_x(x)$ is not primitive recursive: there is no universal primitive recursive function. However, $\psi_x(x)$ for primitive recursive functions is perfectly well defined — there are even recursive functions that compute it and thus recursive functions that can compute the corresponding $h(x)$. The fact that a primitive recursive function for $h(x)$ would lead to contradiction does not stop there being a perfectly concrete Turing machine to compute it.

Similar constructions, showing that $h(x,y)$ or $g(x)$ for Turing machines can be computed by more powerful machines are dealt with routinely in recursion theory (see Rogers [1967] or Odifreddi [1989]) and are not problematic. The fact that hypercomputation is discussed there almost exclusively in terms of oracles and relative computation is of no particular concern: one can simply imagine the set of *o*-machines with a given oracle as a concrete type of hypermachine. In general, the fact that there is no machine of a given type to compute $g(x)$ for that type does not mean that there is any logical contradiction in other machines being able to compute it.

## 2  Issues with Specific Hypermachines

We now turn to some potentially more dangerous claims that are leveled against specific models of hypercomputation. For example, Svozil [1998] claims that accelerating Turing machines can compute *their own* halting function and thus diagonalise out of their set of computable functions. If true, this would certainly require major changes to the model.

An *accelerating Turing machine* (see Copeland [2002]) is a Turing machine with alphabet {1, Blank} whose first square of the tape is a special output square that can be written on once, but never erased. Using some unspecified method, the accelerating Turing machine performs its first operation in half a second, its next in a quarter of a second, then an eighth of a second and so on. Since $1/2 + 1/4 + 1/8 + \ldots = 1$, it can perform an infinity of operations in one second (or some other finite interval). At the end of this time, its output is considered to be whatever is on the special square and this can be checked by an operator.

Such a machine can solve the halting problem for Turing machines by simulating the given Turing machine on its input and marking the special square with a 1 if it halts. An operator can check to see if the given Turing machine halts on its input, by examining the special square after one second. If the square is marked, the answer is yes. If it is unmarked, the answer is no.

It is important to note that an accelerating Turing machine has two relevant timescales: an internal one in which it may take infinitely many timesteps and an external one in which it is always finished by a fixed time. It is with respect to this external timescale that we judge an accelerating Turing machine to compute a function and thus, like primitive recursive functions, accelerating Turing machines compute only total functions and can trivially solve their own halting problems.

However, in an important contrast to primitive recursive functions, accelerating Turing machines are universal. A universal accelerating Turing machine is much like a universal Turing machine, simulating its input machine and marking its special output square if and only if the input machine does. Universality would be enough to allow primitive recursive functions to compute their own $h(x)$, so this seems to put accelerating Turing machines on shaky ground. Can accelerating Turing machines compute their own $h(x)$? Since their output is restricted to either 0 or 1, they obviously cannot, but we can consider a slight modification:

$i(x)$:     $\neg\psi_x(x)$     (i.e. $1 - \psi_x(x)$ )

Given that an accelerating Turing machine can compute $\psi_x(x)$ and that they can also compute negation, what could possibly prevent them from computing $\neg\psi_x(x)$? The answer is that unlike most for models of computation, the functions computed by accelerating Turing machines are not closed under composition. As is quite easily demonstrated (see Ord [2002]), a function is computable by an accelerating Turing machine if and only if it is the characteristic function of a recursively enumerable set. Since this class of functions is not closed under negation, $i(x)$ can remain uncomputable.

A fundamental reason for this lack of closure under composition is that the special output square can only be changed once in a computation and since all non-trivial computations on an accelerating Turing machine involve the potential for changing the output square, they cannot in general be composed. This difficulty can be avoided if the computations run for a finite number of timesteps because normal squares can stand in for the output square, but when a computation needs infinitely many timesteps the special output square must be used and composition is not possible.

There is still a potential problem, however, regarding accelerating Turing machines solving their own *internal halting problem*. Due to the similarity of the two models, the accelerating Turing machine that can solve the halting problem for standard Turing machines can also tell whether an accelerating machine will run for infinitely many internal timesteps. If an accelerating Turing machine could be set up to run for finitely many timesteps if and only if its input machine runs for infinitely many timesteps, then we would have another contradiction. Fortunately however, this is not possible as an accelerating Turing machine cannot in general determine that another accelerating Turing machine runs for infinitely many timesteps without itself taking infinitely many timesteps, and being no longer able to halt within a finite internal time.

Lack of closure under composition makes accelerating Turing machines seem an implausible model of intuitively effective computation or physical computation. This failing can be remedied, however, by embedding them in a Turing complete system. For example, we could treat them as oracles in an *o*-machine, allowing the computation of all functions that are recursive in the halting function. Such a system is then closed under composition and, while it can still solve the halting problem for Turing machines, it has become sufficiently complex to no longer be able to solve its own halting problem and thereby avoid any problems of diagonalisation. Regardless of how we look at accelerating Turing machines — whether in an internal sense, an external sense, or embedded in an *o*-machine — the diagonal method exposes no logical inconsistencies.

Similar to accelerating Turing machines are the *infinite time Turing machines* (ITTMs) of Hamkins and Lewis [1998]. These are like a transfinite extension of accelerating Turing machines, with the configuration of the machine at a successor ordinal timestep being defined from the previous configuration and at a limit ordinal being defined from the set of previous configurations. Among the many abilities of ITTMs is that of deciding the truth of any statement of first order arithmetic (in less than $\omega^2$ timesteps). Cotogno [2003] states that since some arithmetical statements can be built around diagonalisation, we are faced with a dilemma of the type discussed earlier: either use an arbitrary value chosen by convention or face inconsistency. Cotogno then points to the arbitrary convention that there is a bias favouring 1 over 0 when defining the tape configuration at a limit timestep, as the way that ITTMs escape this dilemma.

While this *is* a strictly conventional bias with no intrinsic reason to prefer 1 to 0, it is not related to decisions an ITTM makes regarding the truth of arithmetical expressions. When an ITTM decides such an expression, there is nothing arbitrary about its output. The truth of falsity of the statement is given by the standard definition of truth in arithmetic due to Tarski (see Odifreddi [1989]), which completely specifies this truth-value (although in a non-constructive way). When an ITTM decides an arithmetical expression, it can be seen to conform to Tarski's definition, assigning the one and only value that this definition specifies. Indeed, if we modify the definition of ITTMs so that the bias is for 0 over 1, and then find one of these modified ITTMs that decides arithmetic, it will assign exactly the same truth-values as the unmodified kind.

Not only does an ITTM have no room for arbitrariness when it decides arithmetic, it does not need any. The workings of a given ITTM cannot be captured in first order arithmetic (or any other system that an ITTM is capable of deciding) and thus ITTMs cannot compute their own halting and diagonal functions. Instead, any diagonal function it can compute is one for a different type of machine and thus does not lead to contradiction.

The ITTM is, however, not the most powerful hypermachine that one encounters. Consider, for example, the processor networks of Siegelmann and Sontag [1994], which can compute all functions from N to N. Might these incredibly powerful machines, be able to compute their own $g(x)$? The answer is once again in the negative.

Processor networks are a form of recurrent neural network with real valued weights and natural number input and output. Siegelmann and Sontag showed that these networks can compute all functions from N to N (although less if restricted to polynomial time). The escape from contradiction comes in the fact that there are uncountably many such networks. While one could say that this immediately invalidates the diagonal argument's form, as the set of functions computable by processor networks is uncountable, this would be ignoring a powerful and closely related diagonal argument.

Let us now consider a very general form of the diagonal arguments above. Instead of a sequence of functions, consider a set of functions $\{\psi_i\}$ taking values in $X$ to values in $Y$ where there is some index set $I$ such that $i \in I$. We must also insist that $I \subseteq X$ to guarantee that there is an input code for each computable function.

$j(x)$:  **if** $(x \notin I$ **or** $\psi_x(x)$ diverges$)$ **then** $y_0$ **otherwise** $k(\psi_x(x))$

where:

$y_0$ is some arbitrary member of $Y$. For $g(x)$ we used 0, but in general $Y$ need not include 0.

and:

$k(x)$ is any function where for all $y \in Y$, $k(y) \neq y$. If $k(y)$ diverges it is considered to not equal $y$.

Specifying $j(x)$ obviously involves choosing a particular $y_0$ and $k(x)$, but here it suffices to let these be anything that satisfies the above conditions.

While processor networks can compute all functions from N to N, they cannot compute any $j(x)$

satisfying these conditions. This is because there is no universal processor network — one that takes an arbitrary processor network as input and simulates its behaviour. Such a network cannot exist due to the (trivial) reason that there are only countably many inputs that a processor network can take, but there are uncountably many processor networks. In other words, no appropriate index set *I* is a subset of their input set N and so the diagonal argument does not produce any uncomputable $j(x)$, sparing processor networks from any potential contradictions.

## 3 Conclusions for Hypercomputation

We have explained how some specific hypermachines avoid any contradictions associated with diagonalisation and in doing so have presented a more general diagonalisation function, $j(x)$. Using this approach we can see, for instance, that while processor networks avoided contradiction by having a small input size, one could conceive of logically consistent hypermachines that had domains large enough to encompass their uncountable index set. An example would be a recurrent neural network with real valued input and output. Some hypermachines such as this might then be universal, forcing a more detailed analysis of why they cannot compute their own $j(x)$.

We can go further and present a unified account of such possibilities and what they imply. In this general form, the diagonal argument shows us that no logically consistent machine (hyper or otherwise) can compute its own $j(x)$. It also shows us a route to computing $j(x)$, so we can consider it to describe a set of properties that no logically consistent machine can possess. For if a class of machines computing from *X* to *Y*:

(1) Has an encoding of each machine as an input,
(2) Can determine whether a given input encodes a machine in the class,
(3) Can determine whether a coded machine diverges when applied to its own code,
(4) Can compute the results of applying a coded machine to its own code (given it halts),
(5) Can compute some function that when applied to $y \in Y$ does not return *y*.
(6) Can perform conditional branches, and
(7) Is closed under composition,

then it could compute its own $j(x)$ and we would have a contradiction. Thus, no machine can possess all properties (1) – (7). Of course, this list of mutually inconsistent properties applies not only to machines, but to any paradigm in which we might discuss computation, such as programs or functions.

While possessing all of (1) – (7) allows the computation of a model's own diagonal function, the converse is not always true. For example, a model may lack conditional branching, but may have some other strange primitive operations that allow the computation of $j(x)$. However, these properties are still very useful for showing how some classes of machines avoid contradiction as well how other classes can suffer from it.

Looking at some computational models, we can see:
- Turing machines do not have property (3).
- Primitive recursive functions do not have property (4).
- Total recursive functions either do not have property (2) or do not have (4), depending on the coding used.
- Accelerating Turing machines do not have property (7).
- *o*-machines with accelerating Turing machines as oracles do not have property (3).

- ITTMs do not have property (3).
- Processor networks do not have property (1) which means they also lack (2), (3) and (4).
- Quantum Adiabatic Computers (as presented in Kieu [2003]) do not have property (1) which means they also lack (2), (3) and (4).

In summing up his arguments regarding diagonalisation and hypercomputation, Cotogno [2003] writes:

'The conclusion is therefore negative: one may expect some progress in efficiency, but none in the enlargement of the set of effectively computable functions; this applies to all theories of computation, so long as they share the same Boolean background and are thus subject to the effects of logical incomputability.'

We reach a very different conclusion. Despite claims to the contrary, diagonalisation does not provide much of a threat to the logical consistency of hypercomputation. In particular, there is no inherent problem with hypermachines computing the halting function for Turing machines or the related diagonal function, since this is only impossible for recursive systems. If a type of hypermachine could compute its own the diagonal function $j(x)$, this would indeed lead to contradiction, but on close examination, few if any hypermachines can do so. This is made even clearer by examining the steps in the diagonal method and the properties needed to complete them. This list of seven incompatible properties of machines provides a useful tool for immediately checking whether a machine can compute its diagonal function and demonstrating why or why not.

Chaitin, Gregory J. [1987]: *Algorithmic Information Theory*. Cambridge University Press, Cambridge.

Copeland, B. Jack [1998]: 'Even Turing Machines Can Compute Uncomputable Functions', pp 150-164. In Cristian S. Calude, J. Casti and M. J. Dineen (eds), *Unconventional Models of Computation*. Springer-Verlag, Singapore.

Copeland, B. Jack [2002]: 'Hypercomputation', *Minds and Machines*, **12**, pp. 461-502.

Cotogno, P. [2003]: 'Hypercomputation and the Physical Church-Turing Thesis', *Brit. J. Phil. Sci.* **54**, pp. 181-223.

Kieu, Tien D. [2003]: 'Computing the non-computable', *Contemporary Physics*, **44**, pp. 51-71.

Hamkins, Joel D. and Andy Lewis [1998]: 'Infinite Time Turing Machines'. *Journal of Symbolic Logic*, **65** pp. 567-604.

Péter, Rózsa [1967]: *Recursive Functions*, Academic Press, London.

Odifreddi, Piergiorgio [1989]: *Classical Recursion Theory*. Elsevier, Amsterdam.

Ord, Toby [2002]: 'Hypercomputation: Computing more than the Turing machine'. Technical report, University of Melbourne, Melbourne, Australia. Available at http://arxiv.org/abs/math.LO/0209332.

Rogers, Hartley [1967]: *Theory of Recursive Functions and Effective Computability*. McGraw-Hill.

Siegelmann, Hava T. and Eduardo D. Sontag [1994]: 'Analog Computation via Neural Networks', *Theoretical Computer Science*, **131**, pp. 331-360.

Svozil, Karl. [1998]: 'The Church-Turing Thesis as a Guiding Principle for Physics', pp 371-385. In Cristian S. Calude, J. Casti and M. J. Dineen (eds), *Unconventional Models of Computation*. Springer-Verlag, Singapore.

Turing, Alan M. [1936]: 'On Computable Numbers, with an Application to the Entscheidungsproblem' *Proceedings of the London Mathematical Society*, **42** pp. 230-265.

Toby Ord,
Department of Philosophy, University of Melbourne, Parkville 3010, Australia.
t.ord@pgrad.unimelb.edu.au

Tien D. Kieu,
Centre for Atom Optics and Ultrafast Spectroscopy, Swinburne University of Technology,
Hawthorn 3122, Australia.
kieu@swin.edu.au